\documentclass[11pt,a4paper]{article}
\usepackage[margin=24mm,headheight=14pt]{geometry}
\usepackage{amsmath,amssymb,amsthm,mathtools}
\usepackage[T1]{fontenc}
\usepackage{lmodern,textcomp}
\usepackage{microtype,graphicx,booktabs,longtable,array,calc}
\usepackage{fvextra,seqsplit,xurl}
\usepackage{tikz,pdflscape,caption,pdfpages}
\usepackage{fancyhdr,needspace}
\usepackage[unicode,hidelinks]{hyperref}
\usepackage{bookmark}
\DefineVerbatimEnvironment{verbatim}{Verbatim}{fontsize=\footnotesize,breaklines=true,breakanywhere=true}
\providecommand{\tightlist}{\setlength{\itemsep}{0pt}\setlength{\parskip}{0pt}}
\newcommand{\codewrap}[1]{{\ttfamily\small\seqsplit{#1}}}
\newcommand{\paperheading}[3]{\par\addvspace{8pt}\Needspace{3\baselineskip}\phantomsection\label{result:#2}\noindent\textbf{#1 #2#3.}\enspace}

\hypersetup{pdftitle={Certified exact bounds for adaptive quantitative group testing with two defectives},pdfauthor={Fedor Karpelevitch}}
\begin{document}
\thispagestyle{plain}
\begin{center}{\Large Certified exact bounds for adaptive quantitative group testing with two defectives\par}\vspace{8pt}
Fedor Karpelevitch\\[3pt]
\small Independent researcher\\
\href{mailto:fedor@karpelevitch.net}{\nolinkurl{fedor@karpelevitch.net}}\end{center}
\begin{abstract}

We determine the maximum population size for which exactly two defective items
can always be identified using at most $k$ adaptive quantitative tests, for every
$0\le k\le10$. Each test reports the number of defectives in a chosen subset.
The exact maxima for eight, nine and ten tests are 65, 112 and 192, respectively.
We also classify all pairs of sizes of two disjoint sets, each containing exactly
one defective, that can be resolved within each test budget through eight.
These exact results are supported by explicit strategies and maximality
certificates checked by programs separate from the search implementation.
Additional results include partial nine-test boundaries and an eleven-test
strategy for 328 items; maximality of the latter remains open.

\end{abstract}



\Needspace{11\baselineskip}
\section{Introduction}\label{sec:1}

Suppose that exactly two items in a finite population are defective. A test selects a subset and reports whether it contains zero, one or two defectives; later tests may depend on earlier answers. We seek strategies minimizing the worst-case number of tests. Equivalently, in the graph formulation of Aigner \cite{A86,A88}, the unknown object is an edge and a test reports how many endpoints lie in the selected vertex set. The complete graph \(K_n\) represents an unrestricted population; \(K_{n,m}\) represents two disjoint sets, each known to contain one defective. The two vertex sets of a complete bipartite graph are called its shores.

A depth-\(k\) ternary decision tree has at most \(3^k\) leaves. This information bound is necessary but does not account for the restrictions on tests of defective pairs. Establishing an exact maximum therefore requires both a strategy at the proposed endpoint and an impossibility proof beyond it. We establish exact finite bounds with explicit strategies and exhaustive refutations, and supply a verification procedure that does not share source with the search program. The problem also has a counterfeit-coin interpretation, but the oracle here returns an exact count in one tested subset, not a comparison between two pans of a balance.

\Needspace{7\baselineskip}
\subsection{Related work}\label{sec:1.1}

Aigner \cite{A86} introduced the graph formulation, determined the complete bipartite thresholds with a fixed shore of size two, three or four, and tabulated small complete-graph and bipartite cases. His complete-graph table, repeated in \cite{A88}, gives 37 as the seven-test maximum. Gargano, Montuori, Setaro and Vaccaro \cite{GMSV92} subsequently gave strategies for \(K_{14,9}\) in five tests, \(K_{21,17}\) in six, and \(K_{32,32}\) in seven. Their \(K_{21,17}\) strategy, together with the known six-test complete-graph threshold, already implies a seven-test construction for \(K_{38}\). Likewise, their \(K_{32,32}\) strategy implies an eight-test construction for \(K_{64}\). Our contribution at seven tests is the matching upper bound and the exact correction to 38, not first achievability. At eight tests we establish the exact maximum 65, improving the latter construction.

The relationship to binary-adder coding supplies further constructions. Zhang, Berger and Massey \cite{ZBM87} give finite full-feedback codes corresponding to complete bipartite search strategies. Hwang \cite{H87,H89} surveys the two-defective problem; Hao \cite{H90} and Christen \cite{C94} develop recursive constructions and bounds. Li, Wu and Triesch \cite{LWT18} determine the exact threshold for a fixed shore of size five. These results provide both prior endpoints and checks on the finite frontiers below. Florin, Ho and Jiang \cite{FHJ22} determine the sharp asymptotic rate through the binary adder channel with complete feedback. Our results concern finite maxima and do not replace or improve that asymptotic theorem.

For star forests, Aigner \cite[Propositions\nobreakspace{}3.24-3.25]{A88} gives an explicit canonical strategy and a necessary weak-majorization condition. These are sufficient for the proof rules used here. His conjectured converse is false, as shown separately in \cite{Kar26a}; neither that counterexample nor the small-budget converse proved there is a premise of the finite bounds in this paper.

\Needspace{7\baselineskip}
\subsection{Results and organization}\label{sec:1.2}

Our first results are the exact complete-graph thresholds through ten tests and the complete bipartite frontiers through eight (Theorems\nobreakspace{}\hyperref[result:3.1]{3.1} and \hyperref[result:3.2]{3.2}). Achievability and maximality are certified separately. Section\nobreakspace{}\hyperref[sec:3.3]{3.3} collects the partial nine-test frontier results, including an additional exact small-shore endpoint. Section\nobreakspace{}\hyperref[sec:3.4]{3.4} gives higher-budget bounds and an eleven-test complete-graph construction.

Section\nobreakspace{}\hyperref[sec:2]{2} fixes the model and elementary reductions, including the canonical star strategies. Section\nobreakspace{}\hyperref[sec:3]{3} states the finite results. Section\nobreakspace{}\hyperref[sec:4]{4} describes their verification and distinguishes the evidence for the additional bounds; Section\nobreakspace{}\hyperref[sec:5]{5} gives open questions. Appendix\nobreakspace{}\hyperref[sec:A]{A} proves the certificate inference rules, Appendix\nobreakspace{}\hyperref[sec:B]{B} specifies proof inputs and implementation notation, Appendix\nobreakspace{}\hyperref[sec:C]{C} gives the complete frontier table, and Appendix\nobreakspace{}\hyperref[sec:D]{D} presents an indexed ten-test strategy map.

\Needspace{11\baselineskip}
\section{Model and preliminary lemmas}\label{sec:2}

For a finite simple graph \(G\), a candidate defective pair is an edge \(e\in E(G)\). Testing \(T\subseteq V(G)\) returns \(|e\cap T|\). Let \(M(G)\) be the minimum worst-case number of tests needed to identify \(e\), following the notation of \cite{A88}. Graphs with at most one edge require no tests; an edgeless graph represents an impossible outcome. We ignore isolated vertices. Write \(\mu(G)=|E(G)|\) for the number of hypotheses, also called the mass of the state. Mass counts candidate pairs, not items.

Starting from \(K_n\), an outcome of one defective in a tested subset leaves one defective on each side of the test, giving a complete bipartite graph. Subsequent tests can split this graph into vertex-disjoint complete bipartite components. We therefore consider states of the form \[
G=\bigsqcup_{i=1}^r K_{n_i,m_i},\qquad \mu(G)=\sum_{i=1}^r n_i m_i.
\] Interchanging the shores of a component or permuting components does not change the state. A test choosing \(a_i\) and \(b_i\) vertices from its two shores produces the following graphs \(G_j\) of candidate pairs after outcome \(j\): \[
\begin{aligned}
G_2&=\bigsqcup_i K_{a_i,b_i},\\
G_0&=\bigsqcup_i K_{n_i-a_i,m_i-b_i},\\
G_1&=\bigsqcup_i\left(K_{a_i,m_i-b_i}\sqcup K_{n_i-a_i,b_i}\right).
\end{aligned}
\] Zero-shore components have no edges and are omitted. For \(K_n\), testing \(a\) vertices gives \(K_{n-a}\), \(K_{a,n-a}\) and \(K_a\) for outcomes zero, one and two, respectively. We call outcomes zero and two pure and outcome one mixed.

Define the threshold functions \[
A(k)=\max\{n:M(K_n)\le k\},\qquad
n(k,m)=\max\{n\ge m:M(K_{n,m})\le k\},
\] the second only when the set is nonempty. Here \(k\) always denotes a test budget; \(K_n\) and \(K_{n,m}\) denote graphs, not budgets or threshold functions.

A component with one shore of size one is a star. Star forests provide explicit terminal strategies and necessary bounds for more general states, as explained in Section\nobreakspace{}\hyperref[sec:2.1]{2.1} and \hyperref[sec:A]{A}ppendix A.1. A star forest is \(F(a)=\bigsqcup_i K_{1,a_i}\), represented by its nonincreasing sequence of positive star sizes \(a=(a_1,a_2,\ldots)\). Such a sequence is an integer partition; its entries are also called rows and their sizes widths. Its mass is \(\sum_i a_i\). We write \(a\preceq_w b\) when \(\sum_{i=1}^t a_i\le\sum_{i=1}^t b_i\) for every \(t\), padding with zeros as necessary. This is weak majorization, not coordinatewise comparison. In partitions, \(x^{[r]}\) denotes \(r\) repeated entries of size \(x\). We use Aigner\textquotesingle s notation \(N(k)\) for his canonical star partition, defined below. Repository state names are confined to the implementation concordance in Appendix\nobreakspace{}\hyperref[sec:B.1]{B.1}.

\paperheading{Lemma}{2.1}{ (strategy restriction and pullback)} If a vertex map from a graph \(H\) to a graph \(G\) sends every edge to an edge and is injective on edges, then \(M(H)\le M(G)\).

\begin{proof} Replace each test $T$ in a strategy for $G$ by its preimage in $H$. Each edge gives
the same response as its image at every node. Distinct edges have distinct image edges, so
the strategy distinguishes them. In particular, taking $H$ to be a subgraph proves
Aigner's subgraph monotonicity, $H\subseteq G\Longrightarrow M(H)\le M(G)$
\cite[equation (3.10)]{A88}. The vertex map may identify nonadjacent vertices; injectivity is
required only on candidate edges. \end{proof}

\paperheading{Lemma}{2.2}{ (disjoint unit edges)} Let \(U\) consist of \(u\) vertex-disjoint edges, also disjoint from \(R\). Then \(M(R\sqcup U)\le k\) if and only if \(M(R)\le k\) and \(\mu(R)+u\le3^k\).

\begin{proof} Necessity follows by restriction and the information bound. For sufficiency,
extend a strategy for $R$ to a full ternary tree of depth $k$, padding short paths with
empty tests. Exactly $\mu(R)$ depth-$k$ words are used by its edges. Assign each unit
edge a distinct unused word. At each node include zero, one or both of that edge's
endpoints to realize its next symbol. The disjoint endpoints let these choices coexist
with the original test and with all other unit edges. The mass condition supplies enough
unused words. \end{proof}

\Needspace{7\baselineskip}
\subsection{Canonical star strategies and necessary bounds}\label{sec:2.1}

\(N(k)\) lists the star sizes in Aigner\textquotesingle s canonical construction \cite[display (3.12)]{A88}. It has \(2^k\) entries and mass \(3^k\). Explicitly, put

\[
V_k(j)=\sum_{i=j}^{k}\binom{k}{i}\qquad(0\le j\le k).
\]

For \(k\ge1\), the partition has one entry \(V_k(0)=2^k\) and \(2^{j-1}\) copies of \(V_k(j)\) for each \(1\le j\le k\):

\[
N(k)=\bigl(2^k,V_k(1)^{[1]},V_k(2)^{[2]},\ldots,V_k(k)^{[2^{k-1}]}\bigr).
\]

The base case is \(N(0)=(1)\). For example, \(N(1)=(2,1)\), \(N(2)=(4,3,1,1)\) and \(N(3)=(8,7,4,4,1,1,1,1)\). The equivalent recursive construction explains how the stars can be searched. If \(N(k-1)=(h_1,\ldots,h_s)\), define the zero-padded sequences

\[
\begin{aligned}
L&=(h_1,0,h_2,0,\ldots,h_s,0),\\
B&=(h_1,\ldots,h_s,0,\ldots,0),\\
R&=(0,h_1,0,h_2,\ldots,0,h_s).
\end{aligned}
\]

Each vector has length \(2s\). Then \(N(k)\) is the nonincreasing rearrangement of \(L+B+R\); the explicit formula follows by induction using Pascal\textquotesingle s identity.

The recurrence constructs a strategy for \(F(N(k))\): each row receives mass from the mixed vector and at most one pure vector, and the three children are copies of \(N(k-1)\) up to zero parts and permutations. Induction from \(N(0)=(1)\) proves solvability \cite[Proposition\nobreakspace{}3.24]{A88}. Thus a star forest whose ordered widths fit coordinatewise into distinct \(N(k)\) rows has a strategy by Lemma\nobreakspace{}\hyperref[result:2.1]{2.1}.

Aigner proved the following necessary majorization condition \cite[Proposition\nobreakspace{}3.25]{A88}: \[
M(F(a))\le k\quad\Longrightarrow\quad a\preceq_w N(k)
\] We use this implication only to reject a state that violates a prefix bound. Passing majorization does not itself supply a strategy. Appendix\nobreakspace{}\hyperref[sec:A.1]{A.1} applies this necessary condition to general bipartite states.

\Needspace{11\baselineskip}
\section{Exact finite results}\label{sec:3}

\Needspace{7\baselineskip}
\subsection{Complete graphs}\label{sec:3.1}

\paperheading{Theorem}{3.1}{ (complete-graph maxima; computer-assisted)} The values of \(A(k)\) for \(0\le k\le10\) are exactly those in Table\nobreakspace{}\hyperref[table:1]{1}.

\begin{table}[htbp]
\centering
\caption{Maximum population size for two defectives and at most $k$ tests.}\label{table:1}
\begin{tabular}[t]{r*{11}{r}rrrrrrrrrrrrrrrrrrrrrr}
\toprule
$k$ & 0 & 1 & 2 & 3 & 4 & 5 & 6 & 7 & 8 & 9 & 10 \\ \midrule
$A(k)$ & 2 & 2 & 3 & 5 & 8 & 13 & 22 & 38 & 65 & 112 & 192 \\
\bottomrule
\end{tabular}
\end{table}

\begin{proof} Positive strategies and exhaustive first-test refutations establish the two
directions at each level. The complete bipartite witnesses supply positive first cuts
through nine tests, with direct zero- and one-edge bases at the initial levels; separate
complete-graph trees give the ten-test construction. For the upper bounds, the certificate
lists every first-test size, up to complementation, and refutes at least one child using
an already established smaller-budget claim. Appendix\nobreakspace{}\hyperref[sec:A.3]{A.3} proves this induction, and
Appendix\nobreakspace{}\hyperref[sec:B]{B} identifies the computational inputs. \end{proof}

In particular, 192 items are solvable in ten tests and 193 are not. Once \(A(9)=112\) has been established, the only first tests of \(K_{193}\) whose larger pure child is not already excluded have sizes \(97\le a\le112\). Their mixed children are \(K_{a,193-a}\) with nine tests remaining. The negative certificate refutes all of them. It also derives the preceding complete-graph bound; it does not import that bound as an unchecked solver answer.

Appendix\nobreakspace{}\hyperref[sec:D]{D} gives an indexed map of a ten-test strategy for \(K_{192}\), with a guide to its compact state and test notation.

\Needspace{7\baselineskip}
\subsection{Complete bipartite graphs}\label{sec:3.2}

\paperheading{Theorem}{3.2}{ (complete bipartite frontiers; computer-assisted)} For \(1\le k\le8\), Table\nobreakspace{}\hyperref[table:C.1]{C.1} gives every defined value of \(n(k,m)\). Each entry is an exact maximum. Beyond the last shore size listed for a level, no normalized pair \(n\ge m\) is solvable.

\begin{proof} Every positive boundary has a checked strategy. Its adjacent negative excludes
all greater first-shore sizes. The first excluded diagonal rules out the remaining
normalized pairs. Shore symmetry and Lemma\nobreakspace{}\hyperref[result:2.1]{2.1} therefore classify the entire domain.
Appendices\nobreakspace{}\hyperref[sec:A.1]{A.1} and \hyperref[sec:A.2]{A.2} prove the inference rules used in these certificates. \end{proof}

All eight frontiers are covered by one certificate package: 130 positive boundary strategies and 138 negative-bound derivations, including the excluded diagonals. The full table is in Appendix\nobreakspace{}\hyperref[sec:C]{C} to keep the statements and their proofs together here.

The published fixed-shore formulas provide independent comparisons. Aigner \cite{A86} proves the cases \(m=2,3,4\); the case \(m=1\) is binary search: \[
\begin{aligned}
n(k,1)&=2^k &&(k\ge0),\\
n(k,2)&=2^k-1 &&(k\ge2),\\
n(k,3)&=2^k-k &&(k\ge3),\\
n(k,4)&=2^k-2k+2 &&(k\ge3).
\end{aligned}
\] Writing \(f(k)=2^k-k(k-3)/2-5\), Li, Wu and Triesch \cite{LWT18} prove \[
n(k,5)=
\begin{cases}
f(k),&4\le k\le8,\\
f(k)+1,&9\le k\le10,\\
f(k)+2,&k\ge11.
\end{cases}
\] These domains respect \(n\ge m\): the five-shore formula at three tests gives the reversed pair \(K_{5,3}\), not a value of \(n(3,5)\). Aigner\textquotesingle s small table also gives \(n(4,6)=7\). The finite endpoints already constructed in \cite{ZBM87,GMSV92} remain prior constructions even when the matching maximality certificate is supplied here.

\Needspace{7\baselineskip}
\subsection{Partial nine-test frontier}\label{sec:3.3}

The nine-test frontier is not yet complete. The published formulas in Section\nobreakspace{}\hyperref[sec:3.2]{3.2} determine its first five small-shore entries. The following extends that exact portion.

\paperheading{Proposition}{3.3}{} The nine-test frontier satisfies \[
n(9,6)=473.
\]

\begin{proof} A canonical strategy for $K_{473,6}$ and an exhaustive refutation of
$K_{474,6}$ establish the equality. Appendix\nobreakspace{}\hyperref[sec:A.4]{A.4} explains why the relaxed search
used for this refutation gives an unconditional negative answer. \end{proof}

Further boundary information is supplied by \[
\begin{aligned}
M(K_{112,80})&=9,& M(K_{112,81})&>9,\\
M(K_{95,94})&=9,& M(K_{95,95})&>9.
\end{aligned}
\] The first pair follows from a checked subtree of the ten-test \(K_{192}\) strategy and a root of the complete-graph negative chain. The second pair is supported by retained solver records, not by an independently replayed certificate in the combined package. Appendix\nobreakspace{}\hyperref[sec:B.3]{B.3} specifies their separate evidence. The positive equalities also use the counting bound: each graph has more than \(3^8=6561\) edges, so eight tests cannot suffice.

The retained search snapshot of 26 September 2026 further gives \[
n(9,94)=95,\qquad n(9,93)=96,\qquad n(9,92)\ge97.
\] The first two equalities use completed positive/adjacent-negative pairs in the solver log. No adjacent rejection is yet available for the third entry. These additional results are solver-backed, not independently replayed certificates. The \(112{:}80\)/\(112{:}81\) pair fixes the smaller shore with the larger shore held at 112; it still does not establish \(n(9,80)=112\), since \(K_{113,80}\) remains unresolved. The ninth frontier is incomplete.

Figure\nobreakspace{}\hyperref[figure:1]{1} relates the bipartite frontiers to the complete-graph thresholds. The first-test decomposition gives \[
A(k+1)=\max\{n+m:0\le m\le n\le A(k),\ M(K_{n,m})\le k\}.
\] Indeed, both pure children must be solvable in \(k\) tests, and the mixed child must lie in the \(k\)-test bipartite region; these conditions also suffice (Appendix\nobreakspace{}\hyperref[sec:A.3]{A.3}). For the displayed levels, a maximizing pair can be chosen with \(n=A(k)\). The vertical segment at that cap and the line \(n+m=A(k+1)\) form each orange tooth. The final tooth uses the \(K_{112,80}\) witness and the separately established upper bound \(A(10)=192\), not a complete nine-test frontier.

\Needspace{7\baselineskip}
\subsection{Higher-budget bounds}\label{sec:3.4}

\paperheading{Proposition}{3.4}{} The separate certificates establish \[
n(10,6)=973,\qquad A(11)\ge328.
\]

\begin{proof} An explicit ten-test strategy for $K_{973,6}$ and an exhaustive refutation
of $K_{974,6}$ establish the equality. Every leaf of the positive witness embeds
in Aigner's canonical strategy; Appendix\nobreakspace{}\hyperref[sec:A.4]{A.4} justifies the negative inference.

For the final bound, test 192 of 328 items. The pure outcomes are $K_{192}$ and
$K_{136}$, both solvable in ten further tests by Theorem\nobreakspace{}\hyperref[result:3.1]{3.1} and restriction. The
mixed outcome $K_{192,136}$ has an explicit checked ten-test strategy. Thus all
three children are solvable. \end{proof}

Maximality of 328 in eleven tests is not asserted. In particular, the eleven-test status of \(K_{329}\) and \(K_{330}\) remains open.

\newgeometry{margin=10mm}
\begin{landscape}
\thispagestyle{empty}
\centering
\captionsetup{font=footnotesize,hypcap=false}
\resizebox{0.90\linewidth}{!}{\input{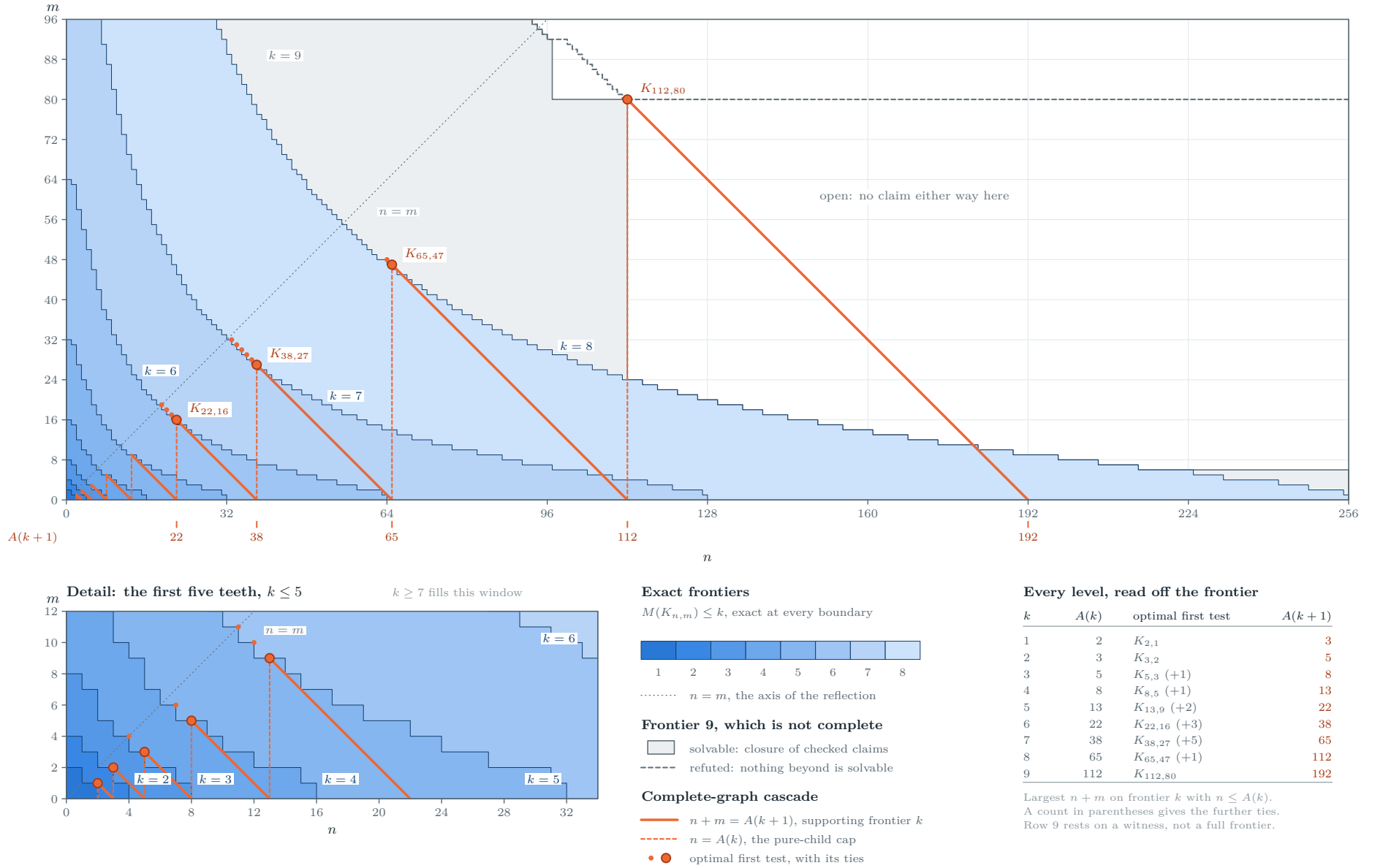}}\par
\captionof{figure}{Complete bipartite frontiers and complete-graph thresholds. Blue regions show solvability through eight tests, with exact boundaries. Shore symmetry is shown on both sides of the dotted diagonal. The grey nine-test region shows established lower bounds; its dashed exclusion boundary leaves an unresolved band. This envelope includes the solver-backed partial frontier described above. Orange teeth show the first-test relation to \(A(k)\); large dots select maximizing cuts and smaller dots show additional known ties. The inset expands the first five teeth. The plot is cropped; the full exact frontiers are in Table\nobreakspace{}\hyperref[table:C.1]{C.1}.}\label{figure:1}
\par\vspace{6pt}\makebox[\linewidth]{\small\thepage}
\end{landscape}
\restoregeometry

\Needspace{11\baselineskip}
\section{Certificates and computational verification}\label{sec:4}

Achievability and maximality have different proof objects. A positive tree exhibits tests that distinguish every candidate pair. A negative certificate supports a separate check that every possible first test has at least one outcome that cannot be resolved with the remaining budget. Those smaller-budget impossibility claims are checked in turn, down to elementary bounds. This recursion certifies exhaustive coverage without trusting the search program\textquotesingle s choice of tests.

A positive tree records a test at every internal node. The checker reconstructs all three children and accepts only proved terminal cases, including coordinatewise embedding into Aigner\textquotesingle s explicit \(N(k)\) strategy. Arbitrary weak majorization is not accepted as a general positive terminal rule. Lemmas\nobreakspace{}\hyperref[result:2.1]{2.1} and \hyperref[result:2.2]{2.2} and Appendix A justify the reductions.

The bipartite frontier certificate has 130 positive trees and 138 negative-bound derivations. Fifty-six negative bounds are roots in a chain containing 5,819,115 claim occurrences; the other 82 follow from information capacity, the star bound, or a cited negative subgraph. The ten-test complete-graph refutation uses a second chain with 2,846,568 occurrences. An explicit complete-graph induction accounts for all 242 complementary first-test classes across eleven negative claims, including the nine-test refutation needed in the final ten-test reduction.

The independent checker is implemented from the mathematical specification and shares no source with the search program. Complete replay reconstructs the split space and requires every assigned claim and requested root to be verified. Full state mass is checked before unit edges are removed. Interrupted runs, unsuccessful workers and structure-only checks are not accepted as complete proofs.

Both negative chains have full semantic replays using the verifier version pinned in the certificate package identified in Appendix\nobreakspace{}\hyperref[sec:B.2]{B.2}, in addition to the positive-witness, reduction and structural checks. The separate bounds in Propositions\nobreakspace{}\hyperref[result:3.3]{3.3} and \hyperref[result:3.4]{3.4} have distinct computational inputs and validation records. The solver-backed 95-shore pair in Section\nobreakspace{}\hyperref[sec:3.3]{3.3} is also separate from the combined certificate. Appendix\nobreakspace{}\hyperref[sec:B]{B} specifies these scopes and reproduction commands.

Table\nobreakspace{}\hyperref[table:2]{2} summarizes the proof objects and completed checks. Here "independent" means a checker that does not share source with the search program; it does not mean a formal proof-assistant verification or an independent human review.

\Needspace{28\baselineskip}

\phantomsection\label{table:2} \emph{Table\nobreakspace{}\hyperref[table:2]{2}. Proof and verification status of the finite results.}

{
\begin{longtable}[]{@{}
  >{\raggedright\arraybackslash}p{(\linewidth - 6\tabcolsep) * \real{0.2500}}
  >{\raggedright\arraybackslash}p{(\linewidth - 6\tabcolsep) * \real{0.2500}}
  >{\raggedright\arraybackslash}p{(\linewidth - 6\tabcolsep) * \real{0.2500}}
  >{\raggedright\arraybackslash}p{(\linewidth - 6\tabcolsep) * \real{0.2500}}@{}}
\toprule\noalign{}
\begin{minipage}[b]{\linewidth}\raggedright
Result
\end{minipage} & \begin{minipage}[b]{\linewidth}\raggedright
Achievability
\end{minipage} & \begin{minipage}[b]{\linewidth}\raggedright
Impossibility or maximality
\end{minipage} & \begin{minipage}[b]{\linewidth}\raggedright
Completed validation
\end{minipage} \\
\midrule\noalign{}
\endhead
\bottomrule\noalign{}
\endlastfoot
Theorems\nobreakspace{}\hyperref[result:3.1]{3.1}-3.2 & Explicit trees & Complete negative chains and boundary reductions & Independent positive checks and full negative replay \\
\(K_{112,80}\) / \(K_{112,81}\) bracket & Checked subtree & Complete-graph chain root & Covered by the independent checks \\
Other partial ninth-frontier claims in Section\nobreakspace{}\hyperref[sec:3.3]{3.3} & Retained constructions & Adjacent rejections where stated & Solver records; no independent negative replay \\
\(n(9,6)\) and \(n(10,6)\) & Canonical trees & Exhaustive relaxed-search rejections & Independent positive checks; negative searches reproduced with the same implementation \\
\(A(11)\ge328\) & Explicit first cut and checked trees & No maximality claim & Independent positive checks \\
\end{longtable}
}

Full negative replays with this verifier took approximately one hour on an AWS c8a.16xlarge using 64 workers for the complete-graph chain and five hours on an AWS c8a.24xlarge using 96 workers for the bipartite chain. Appendix\nobreakspace{}\hyperref[sec:B.2]{B.2} gives their resource measurements and scope, separately from the preliminary checks.

\Needspace{11\baselineskip}
\section{Discussion and open questions}\label{sec:5}

The next finite problems are to complete the nine-test bipartite frontier and determine the eleven-test complete-graph maximum. The present theorems do not require either problem to be settled.

\paperheading{Conjecture}{5.1}{ (antidiagonal transfer)} For \(n\ge m>1\) and \(k\ge0\), \[
M(K_{n,m})\le k\quad\Longrightarrow\quad M(K_{n+1,m-1})\le k.
\] This implication is not used in the proofs or certificates. It reduces the number of edges by \(n-m+1\), but changes both shores and supplies no subgraph embedding of the kind needed by Lemma\nobreakspace{}\hyperref[result:2.1]{2.1}.

Finite patterns alone do not establish general formulas or asymptotic bounds. For example, the expression \(2^k-k(k-1)/2-3\) for \(n(k,6)\) agrees with the exact values for \(4\le k\le9\), but predicts 976 at ten tests, contradicting Proposition\nobreakspace{}\hyperref[result:3.4]{3.4}. Any new asymptotic claim would also have to be compared with the sharp rate already established in \cite{FHJ22}.

First-test selection heuristics can favor small mixed outcomes after qualifying the pure outcomes. Such preferences are not completeness principles: a smaller number of candidate pairs need not make a graph easier to resolve. We use these objectives to order candidate constructions, while exhaustive search and certificate verification remain separate from the ordering choices.

\section*{Data and code availability}

Reproduction materials are available at \href{https://doi.org/10.5281/zenodo.23006289}{10.5281/zenodo.23006289}. Appendix\nobreakspace{}\hyperref[sec:B]{B} identifies the certificates, verifier sources and execution records. The combined finite-certificate package covers Theorems\nobreakspace{}\hyperref[result:3.1]{3.1}\textendash{}\hyperref[result:3.2]{3.2}; the additional results in Sections\nobreakspace{}\hyperref[sec:3.3]{3.3}\textendash{}\hyperref[sec:3.4]{3.4} require the separate inputs specified there.

\section*{Acknowledgments}

The author thanks Konstantin Knop and Igor Krivokon for introducing the problem and for their advice and consultations over the years.

\section*{Funding and competing interests}

This research received no external funding. The author declares no competing interests.

\section*{Use of generative AI}

Generative-AI assistants, primarily OpenAI Codex and, to a lesser extent, Anthropic Claude, were used for literature search, development and exposition of mathematical arguments, software development and debugging, computational experimentation, and manuscript and figure preparation. The author directed this work and assumes responsibility for all content, including the mathematical claims, references, software and reported computations. AI output is not itself evidence for a mathematical claim; the arguments and computational verification procedures, with their scopes and limitations, are described in the paper and appendices.

\clearpage
\appendix
The appendices use the graph model and notation of the main text. Appendix\nobreakspace{}\hyperref[sec:A]{A} proves the inference rules; Appendix\nobreakspace{}\hyperref[sec:B]{B} identifies their computational inputs and reproduction scope; Appendix\nobreakspace{}\hyperref[sec:C]{C} lists the complete certified frontiers.

\Needspace{11\baselineskip}
\section{Soundness of the finite certificates}\label{sec:A}

\Needspace{7\baselineskip}
\subsection{Positive trees and star refutations}\label{sec:A.1}

\paperheading{Lemma}{A.1}{ (full star expansion)} For \(G=\bigsqcup_{i=1}^r K_{n_i,m_i}\), orient \(n_i\ge m_i\) and let \(\Phi(G)\) be the sorted sequence containing \(m_i\) copies of \(n_i\) for each part. If \(G\) is solvable in \(k\), then \(\Phi(G)\preceq_w  N(k)\).

\begin{proof} Replace each $K_{n,m}$ by $m$ disjoint stars of width $n$. Map each star
centre to its vertex on the $m$ shore and its $n$ leaves to the corresponding vertices
on the $n$ shore. The induced edge map is bijective. Lemma\nobreakspace{}\hyperref[result:2.1]{2.1} pulls back a strategy
for the original state; star-forest necessity then gives the inequality.

This is a vertex-splitting relaxation, not a subgraph inclusion: leaves belonging
to different stars may map to the same original vertex. Distinct expanded edges
still map to distinct original edges. Thus a strategy for $G$ induces one for the
expanded forest, but not necessarily conversely. This direction is sufficient:
an unsolvable expanded forest excludes a solution of $G$. \end{proof}

This preserves the full mass, unlike keeping just one star per component. A failed prefix inequality is a negative certificate. Passing the test supplies no positive certificate. After \(2^k\) entries its right side is the constant \(3^k\), so the remaining comparison is already covered by the information bound.

A positive tree is verified by induction on its remaining test budget. At a split node the selected sizes lie within their parent shores, and all three children must equal the graphs computed by the split formulas in Section\nobreakspace{}\hyperref[sec:2]{2} of the manuscript. An internal edge consumes one test. A terminal may have at most one candidate edge or may be a star forest with its sorted widths fitting coordinatewise into distinct \(N(k)\) rows. The explicit canonical strategy and Lemma\nobreakspace{}\hyperref[result:2.1]{2.1} justify the latter. A forest merely weakly majorized by \(N(k)\) is not an accepted general positive terminal.

\Needspace{7\baselineskip}
\subsection{Negative bipartite chains}\label{sec:A.2}

A negative claim at level \(k\) is discharged directly by information capacity or by Lemma\nobreakspace{}\hyperref[result:A.1]{A.1}, or by showing that every legal first test has an unsolvable child at \(k-1\). A child can be refuted by either direct condition or by a proved negative subgraph. One sufficient subgraph relation injects distinct stored components into distinct child components, with both shores no larger after orientation. Multiplicities and skipped child components must be handled explicitly.

The unit rule of Lemma\nobreakspace{}\hyperref[result:2.2]{2.2} requires full mass to be checked before stripping unit edges. In particular a negative support state above information capacity cannot contribute its stripped core as a negative fact: that core might be solvable. An admitted negative support state within capacity has a negative core, by the contrapositive of the unit rule.

For exhaustive coverage, each component \(K_{n,m}\) has exactly the legal choices \(0\le a\le n\), \(0\le b\le m\). Square-shore interchange may normalize its option. Options on equal components may then be sorted nondecreasingly, retaining repetitions. Global test complementation interchanges pure outcomes and preserves the mixed graph. Apply it to the entire option vector, normalize squares and equal runs again, and retain the lexicographically smaller vector, including equality. This is an involution on the already normalized representatives, so no combined symmetry orbit disappears.

A partial child\textquotesingle s negative subgraph persists in every extension of the chosen prefix, justifying early rejection of that prefix. A completed split must have at least one discharged child. Hints may change enumeration order but cannot remove a remaining split. Induction on the level then proves every claim when the support at each level is exactly the already proved set below it and the lowest support is empty. Acceptance requires that every assigned claim was processed, every worker succeeded, and every requested root is present. Structural consistency alone does not establish this induction\textquotesingle s semantic premise.

For a complete frontier, positive boundaries prove every smaller first shore by restriction. The adjacent negatives rule out larger first shores. If \(d\) is the last feasible diagonal, the excluded state \(K_{d+1,d+1}\) rules out all remaining normalized pairs. Thus the positive/negative boundaries plus one excluded diagonal classify the entire domain.

\Needspace{7\baselineskip}
\subsection{Complete-graph induction}\label{sec:A.3}

Every first test on \(K_{N}\) is, up to relabeling and complementation, described by \(\lceil N/2\rceil\le a\le N\), with \(b=N-a\). The children are \(K_{a}\), \(K_{a,b}\) and \(K_{b}\) at level \(k-1\). The endpoint \(b=0\) remains a legal, uninformative test and is included. If \(K_{t}\) is already refuted at \(k-1\) and \(t\le a\), its negative subgraph refutes the large pure child. Otherwise the mixed child can be refuted by information capacity or by a certified negative bipartite subgraph. Covering all first-test sizes proves \(K_{N}\) negative at \(k\).

The combined certificate starts with the information refutation \(M(K_3)>0\) and includes every required level through ten. Each pure citation points to an already established claim at the preceding level. In particular it proves \(M(K_{113})>9\) before using that state to refute the large pure child of \(K_{193}\) at depth ten for \(a\ge113\). For the remaining classes \(97\le a\le 112\), the second bipartite chain refutes all sixteen mixed roots \(K_{a,193-a}\) in nine tests. These two ranges cover all first tests without an external complete-graph premise.

For achievability, if \(a,b\le A(k-1)\) and \(K_{a,b}\) is solvable at \(k-1\), a first test on \(a\) of \(a+b\) items proves \(A(k)\ge a+b\). For any positive integers \(a,b\) this follows by the same three-child decomposition. More generally, allowing empty first shores and treating an empty mixed graph as solved gives the exact recursion

\[
A(k)=\max\{a+b:0\le a,b\le A(k-1),\ M(K_{a,b})\le k-1\}.
\]

Necessity follows from the first test of a solution; sufficiency is the construction just given. The frontier positive trees supply the lower constructions through nine, with direct zero/one-edge bases at the initial levels; the separate \(K_{192}\) trees give ten. Choosing \(a=192,b=136\) and the recovered ten-test bipartite tree proves \(A(11)\ge 328\). No refutation at eleven tests is included.

\Needspace{7\baselineskip}
\subsection{The separate small-side bounds}\label{sec:A.4}

The retained \(m=6\) computations use a different exhaustive recurrence. At a star forest passing majorization it may stop positively; otherwise it explores all legal rectangle cuts to the available depth. At full depth this is a relaxation of solvability: any real strategy passes it by following its actual cuts and the necessary terminal conditions. Therefore a completed \texttt{NO} refutes solvability, while a \texttt{YES} requires a separately valid positive tree.

This explains why the retained nine-test rejection of \(K_{474,6}\) and the canonical 473 tree prove \(n(9,6)=473\). At ten tests, the rejection of \(K_{974,6}\) combines with a separately checked \(K_{973,6}\) witness to prove \(n(10,6)=973\). The latter expands every nonembedded leaf of the original relaxed tree into an explicit continuation; all its remaining terminals embed in distinct canonical slots. The small-side rejections have separate source and execution records; the combined frontier/complete-graph checker does not certify them merely by completing its own chains.

\Needspace{11\baselineskip}
\section{Computational premises and reproduction}\label{sec:B}

\Needspace{7\baselineskip}
\subsection{Implementation concordance}\label{sec:B.1}

The article uses graph notation throughout. The retained source code and certificate files use older names; they are not additional mathematical objects.

{
\begin{longtable}[]{@{}
  >{\raggedright\arraybackslash}p{(\linewidth - 4\tabcolsep) * \real{0.3333}}
  >{\raggedright\arraybackslash}p{(\linewidth - 4\tabcolsep) * \real{0.3333}}
  >{\raggedright\arraybackslash}p{(\linewidth - 4\tabcolsep) * \real{0.3333}}@{}}
\toprule\noalign{}
\begin{minipage}[b]{\linewidth}\raggedright
Mathematical notation
\end{minipage} & \begin{minipage}[b]{\linewidth}\raggedright
Repository notation
\end{minipage} & \begin{minipage}[b]{\linewidth}\raggedright
Meaning
\end{minipage} \\
\midrule\noalign{}
\endhead
\bottomrule\noalign{}
\endlastfoot
\(K_n\) & \texttt{Sa(n)} & All unordered pairs from one population \\
\(K_{n,m}\) & \texttt{Sb(n:m)} & One defective in each of two disjoint sets \\
\(\bigsqcup_i K_{n_i,m_i}\) & \texttt{Sb(n1:m1,...)} & Disjoint union of bipartite candidate graphs \\
\(F(a)=\bigsqcup_i K_{1,a_i}\) & singleton \texttt{Sb(a1:1,...)} state & A star forest \\
\(N(k)\) & \texttt{G\_k} & Aigner\textquotesingle s canonical partition \\
\(k\) & \texttt{K} or \texttt{k} & Remaining test budget \\
\(M(G)\le k\) & \texttt{can\ solve\ ...\ in\ k} & A positive feasibility claim, subject to its evidence \\
\(A(k)\) & \texttt{pareto\_sa.csv}: \texttt{k,n,bound} & The \texttt{n} value is exact only when \texttt{bound=max} \\
\(n(k,m)\) & \texttt{pareto\_sb.csv}: \texttt{k,m,n1,bound} & The \texttt{n1} value is exact only when \texttt{bound=max} \\
\end{longtable}
}

A file label such as \texttt{Sa(192)@10} specifies a graph and test budget, not a maximum. A negative claim corresponds to \(M(G)>k\); a deadline or an inconclusive search does not establish it. Filenames and literal checker output below retain their original spelling so that the instructions apply to the archived artifacts.

\Needspace{7\baselineskip}
\subsection{Complete finite frontier and complete-graph certificate}\label{sec:B.2}

The immutable release \texttt{complete-certificates-2026-09-13} carries the complete \(k\le 8\) frontiers and the complete-graph induction through ten. Its ZIP SHA-256 is \codewrap{70185a7f1737bc9317dbaddb2c444ea690326d37f07412838cd106bbdd50e82d}. From the extracted directory, the complete replay command is

\begin{verbatim}
python3 checkers/check_sa_certificate.py .
\end{verbatim}

Require exit zero and \texttt{COMPLETE\_RESULT\ status=verified}. The runner checks manifest hashes, positive witnesses, all reduction citations and both full negative chains. The \texttt{-\/-structure-only} alternative explicitly leaves negative semantics unchecked. The assembled package was independently downloaded and hash-checked. A fresh full replay of its Pareto component completed on 2026-09-14: all 5,819,115 negative claims, 130 positive trees and 138 frontier bounds passed. Shared proof and Rust source hashes match this combined package. The separate complete-graph chain\textquotesingle s full post-repair replay completed on 2026-09-27: all 2,846,568 claims passed with exact roots, zero gaps and successful exit. Its executable is byte-identical to the Pareto replay\textquotesingle s verifier. The precise provenance is in the complete-graph certificate record (\texttt{evidence/sa\_refutations\_2026-09-13.md}) and the full Sa (\texttt{evidence/sa\_certificate\_aws\_2026-09-27.md}) and Pareto (\codewrap{evidence/pareto\_certificate\_aws\_2026-09-13.md}) replay records.

\textbf{Time and hardware.} Replay requires a 64-bit CPU, Python 3.9 or later, Rust/Cargo 1.85 or later, and a host C linker; no GPU or solver service is needed. With the toolchain installed, the build is offline and the extracted ZIP requires no zstd. The archive is approximately 67 MB, expanding to 409 MB before compilation; additional disk space is needed for the toolchain, build products and logs.

The quick \texttt{-\/-structure-only} check, which also checks positive trees and reduction citations, took approximately 14\textendash{}20 wall seconds on the local Apple M4 Pro, including the offline build, with less than 1 GB sampled process footprint. This is not the cost of verifying the negative proofs. The fresh full Pareto package replay took approximately five hours on an AWS c8a.24xlarge (AMD EPYC 9R45, 96 CPU cores, up to 96 workers), with about 7 GiB peak RSS. This includes all positive/coverage checks and full negative replay, excluding toolchain installation, compilation, selftest and archival. The separate complete-graph chain\textquotesingle s full repaired-verifier replay took approximately one hour on an AWS c8a.16xlarge (AMD EPYC 9R45, 64 cores/workers; about 3 GiB peak RSS), excluding setup, startup checks and archival. These observed times are not extrapolations, and their sum is not a measured combined runtime.

The current combined full runtime and peak memory remain unmeasured. Plan for at least 16 GB of available RAM as initial headroom, not a demonstrated minimum or a guarantee of completion. The two chains run sequentially; \texttt{-\/-threads\ N} controls concurrency, while \texttt{-\/-max-index-nodes\ N} limits the index but not total process memory. The exact measurements and their sources are recorded in resource planning (\texttt{evidence/sa\_refutations\_2026-09-13.md}, resource-planning-for-readers). These timings exclude the separate small-shore computations.

\Needspace{7\baselineskip}
\subsection{Separate proof inputs}\label{sec:B.3}

The following claims use inputs separate from the combined release, with their own source identities and retained completion records:

\begin{itemize}
\tightlist
\item
  The canonical \texttt{witnesses/canon\_473\_6\_at9.tree} and the retained 474 rejection in \texttt{small-m-frontier-2026-08-15:m6\_k9.log}, summarized in nine-test small-shore record (\texttt{evidence/sb\_m6\_k9\_frontier.txt}).
\item
  The positive \texttt{witnesses/sb-973-6-k10.tree}, documented in the ten-test construction record (\texttt{evidence/sb973\_witness\_2026-09-25.md}), and the retained 974 rejection, source \texttt{tools/search\_singletonization.cpp} and validation record ten-test small-shore record (\texttt{evidence/sb\_m6\_k10\_frontier.txt}), for the separate ten-test maximum. The historical negative record predates the later complete provenance banners. Both small-shore rejections were freshly reproduced for this supplement with complete provenance; their search counts agree with the retained records.
\item
  The recovered \texttt{witnesses/sb-192-136-k10.tree} and a ten-test \(K_{192}\) tree, supporting the eleven-test construction. Its provenance and direct checker output are in 328-item construction record (\texttt{evidence/sa328\_witness\_2026-09-13.md}).
\item
  The partial nine-test pairs in Section\nobreakspace{}\hyperref[sec:3.3]{3.3}, documented in the boundary evidence record (\texttt{evidence/partial\_pareto9\_2026-09-13.md}). The \(K_{112,80}\) subtree is already in the checked \(K_{192}\) witness, and the \(K_{112,81}\) refutation is in the combined negative chain. The \(K_{95,94}\) positive is retained in the archived bootstrap log. The 26 September snapshot (\texttt{evidence/pareto9\_snapshot\_2026-09-25.md}) retains the full later v4 log, including the \(K_{95,95}\), \(K_{96,94}\) and \(K_{97,93}\) rejections and the \(K_{96,93}\) and \(K_{97,92}\) constructions. These remain solver-backed records, not independently replayed certificates in the combined ZIP.
\end{itemize}

The standalone positive check from a repository checkout is

\begin{verbatim}
python3 tools/check_witness.py witnesses/canon_473_6_at9.tree \
  witnesses/sb-973-6-k10.tree \
  witnesses/sb-192-136-k10.tree witnesses/sa192_k10_b.tree
\end{verbatim}

The small-shore records and trees, the recovered eleven-test mixed tree and the solver-backed ninth-frontier records are not contained in the immutable combined ZIP. Its complete runner therefore does not verify those claims. They are collected separately in this paper\textquotesingle s versioned supplement, alongside the unchanged combined package. From that supplement, \texttt{python3\ verify.py\ quick\ -\/-work\ ../checks} checks integrity, positive witnesses and finite proof structure, not the large negative semantics. The separate modes \texttt{small-shore} and \texttt{finite} reproduce the corresponding full computations; none independently certifies the partial ninth frontier. The supplement\textquotesingle s README specifies prerequisites, completion conditions and retained-versus-fresh validation.

\Needspace{11\baselineskip}
\section{Complete bipartite frontiers through eight tests}\label{sec:C}

\begin{table}[htbp]
\centering
\renewcommand{\thetable}{C.1}
\caption{Exact values of \(n(k,m)\) for \(1\le k\le8\), generated from the certified frontier data. A dash means that no pair with \(n\ge m\) is solvable at that level, not that the entry is unknown. Each numerical entry is a maximum.}\label{table:C.1}
{\fontsize{8.5}{11}\selectfont\setlength{\tabcolsep}{2.8pt}
\begin{tabular}[t]{r*{8}{r}rrrrrrrrrrrrrrrr}
\toprule
$m$ / $k$ & 1 & 2 & 3 & 4 & 5 & 6 & 7 & 8 \\ \midrule
1 & 2 & 4 & 8 & 16 & 32 & 64 & 128 & 256 \\
2 & -- & 3 & 7 & 15 & 31 & 63 & 127 & 255 \\
3 & -- & -- & 5 & 12 & 27 & 58 & 121 & 248 \\
4 & -- & -- & 4 & 10 & 24 & 54 & 116 & 242 \\
5 & -- & -- & -- & 9 & 22 & 50 & 109 & 231 \\
6 & -- & -- & -- & 7 & 19 & 46 & 104 & 225 \\
7 & -- & -- & -- & -- & 17 & 42 & 97 & 214 \\
8 & -- & -- & -- & -- & 15 & 38 & 91 & 206 \\
9 & -- & -- & -- & -- & 14 & 36 & 87 & 198 \\
10 & -- & -- & -- & -- & 12 & 33 & 82 & 189 \\
11 & -- & -- & -- & -- & 11 & 31 & 77 & 182 \\
12 & -- & -- & -- & -- & -- & 29 & 73 & 174 \\
13 & -- & -- & -- & -- & -- & 27 & 69 & 168 \\
14 & -- & -- & -- & -- & -- & 25 & 66 & 161 \\
15 & -- & -- & -- & -- & -- & 24 & 63 & 155 \\
16 & -- & -- & -- & -- & -- & 22 & 60 & 150 \\
17 & -- & -- & -- & -- & -- & 21 & 58 & 144 \\
18 & -- & -- & -- & -- & -- & 20 & 55 & 139 \\
19 & -- & -- & -- & -- & -- & 19 & 53 & 135 \\
20 & -- & -- & -- & -- & -- & -- & 51 & 130 \\
21 & -- & -- & -- & -- & -- & -- & 49 & 126 \\
22 & -- & -- & -- & -- & -- & -- & 47 & 122 \\
23 & -- & -- & -- & -- & -- & -- & 45 & 118 \\
24 & -- & -- & -- & -- & -- & -- & 43 & 115 \\
25 & -- & -- & -- & -- & -- & -- & 41 & 111 \\
26 & -- & -- & -- & -- & -- & -- & 40 & 108 \\
27 & -- & -- & -- & -- & -- & -- & 38 & 105 \\
28 & -- & -- & -- & -- & -- & -- & 37 & 102 \\
\bottomrule
\end{tabular}
\hspace{10pt}
\begin{tabular}[t]{r*{8}{r}rrrrrrrrrrrrrrrr}
\toprule
$m$ / $k$ & 1 & 2 & 3 & 4 & 5 & 6 & 7 & 8 \\ \midrule
29 & -- & -- & -- & -- & -- & -- & 36 & 100 \\
30 & -- & -- & -- & -- & -- & -- & 35 & 97 \\
31 & -- & -- & -- & -- & -- & -- & 34 & 94 \\
32 & -- & -- & -- & -- & -- & -- & 33 & 92 \\
33 & -- & -- & -- & -- & -- & -- & -- & 89 \\
34 & -- & -- & -- & -- & -- & -- & -- & 87 \\
35 & -- & -- & -- & -- & -- & -- & -- & 85 \\
36 & -- & -- & -- & -- & -- & -- & -- & 83 \\
37 & -- & -- & -- & -- & -- & -- & -- & 81 \\
38 & -- & -- & -- & -- & -- & -- & -- & 79 \\
39 & -- & -- & -- & -- & -- & -- & -- & 77 \\
40 & -- & -- & -- & -- & -- & -- & -- & 76 \\
41 & -- & -- & -- & -- & -- & -- & -- & 74 \\
42 & -- & -- & -- & -- & -- & -- & -- & 72 \\
43 & -- & -- & -- & -- & -- & -- & -- & 71 \\
44 & -- & -- & -- & -- & -- & -- & -- & 69 \\
45 & -- & -- & -- & -- & -- & -- & -- & 68 \\
46 & -- & -- & -- & -- & -- & -- & -- & 66 \\
47 & -- & -- & -- & -- & -- & -- & -- & 65 \\
48 & -- & -- & -- & -- & -- & -- & -- & 64 \\
49 & -- & -- & -- & -- & -- & -- & -- & 62 \\
50 & -- & -- & -- & -- & -- & -- & -- & 61 \\
51 & -- & -- & -- & -- & -- & -- & -- & 60 \\
52 & -- & -- & -- & -- & -- & -- & -- & 59 \\
53 & -- & -- & -- & -- & -- & -- & -- & 58 \\
54 & -- & -- & -- & -- & -- & -- & -- & 57 \\
55 & -- & -- & -- & -- & -- & -- & -- & 56 \\
\bottomrule
\end{tabular}
}
\end{table}

The machine-readable values and per-entry sources are in \texttt{data/pareto\_sb.csv}. Only the exact levels covered by Theorem\nobreakspace{}\hyperref[result:3.2]{3.2} are reproduced here; higher-level exploratory records are not part of this table.

\clearpage
\Needspace{11\baselineskip}
\section{An indexed ten-test strategy for 192 items}\label{sec:D}

Figure\nobreakspace{}\hyperref[figure:2]{2} presents the checked strategy for \(K_{192}\) used in Section\nobreakspace{}\hyperref[sec:3.1]{3.1}. It is a directed acyclic graph: several outcomes may reuse one continuation. The map retains all 148 test nodes with at least two tests remaining; the single one-test node and immediate terminal continuations are omitted from the drawing, not from the source witness. The full witness is \texttt{witnesses/sa192\_k10\_b.tree}.

A round-ended node labelled \texttt{n{[}t{]}} denotes \(K_n\) and a test of \(t\) vertices. Its next row gives the continuations for outcomes one and two. Outcome zero is implicit: the source\textquotesingle s pure children are \(K_{n-t}\) and \(K_t\), and \(n-t\le t\), so restriction supplies the smaller pure strategy. Plain numeric references identify these complete-graph nodes.

A rectangular node has three lines. The first gives its identifier followed by continuations for outcomes zero, one and two, in that order. The second, \texttt{(n1:m1,n2:m2,...)}, represents \(\bigsqcup_i K_{n_i,m_i}\). The third, \texttt{{[}a1:b1,a2:b2,...{]}}, gives the selected shore sizes in the same component order. The leading digit of each rectangular-node identifier is the remaining test budget; the letters distinguish nodes within that level.

A plain continuation is an exact outcome state up to graph symmetries. A tilde marks a checked reuse of another strategy, including the unit-edge reduction of Lemma\nobreakspace{}\hyperref[result:2.2]{2.2} and its mass check; it need not mean containment of the entire printed state. A dash denotes an omitted one-test or immediate terminal continuation, not necessarily a state already resolved. Only selected links are drawn to preserve legibility; other destinations remain explicit in the node labels. Solid curves show exact continuations and dotted curves show reuse. Blue, amber and green distinguish outcomes zero, one and two; crossings are not junctions. Every displayed nonroot node has a drawn incoming link, and every displayed nonterminal has a drawn outgoing link.

The following landscape page is intended as a zoomable reference map. It illustrates the construction, but the checked witness, rather than the drawing, is the proof object.

\clearpage
\includepdf[pages=1,fitpaper=true,pagecommand={\thispagestyle{empty}\refstepcounter{figure}\label{figure:2}},addtotoc={1,subsection,2,{Ten-test strategy map},strategy-map}]{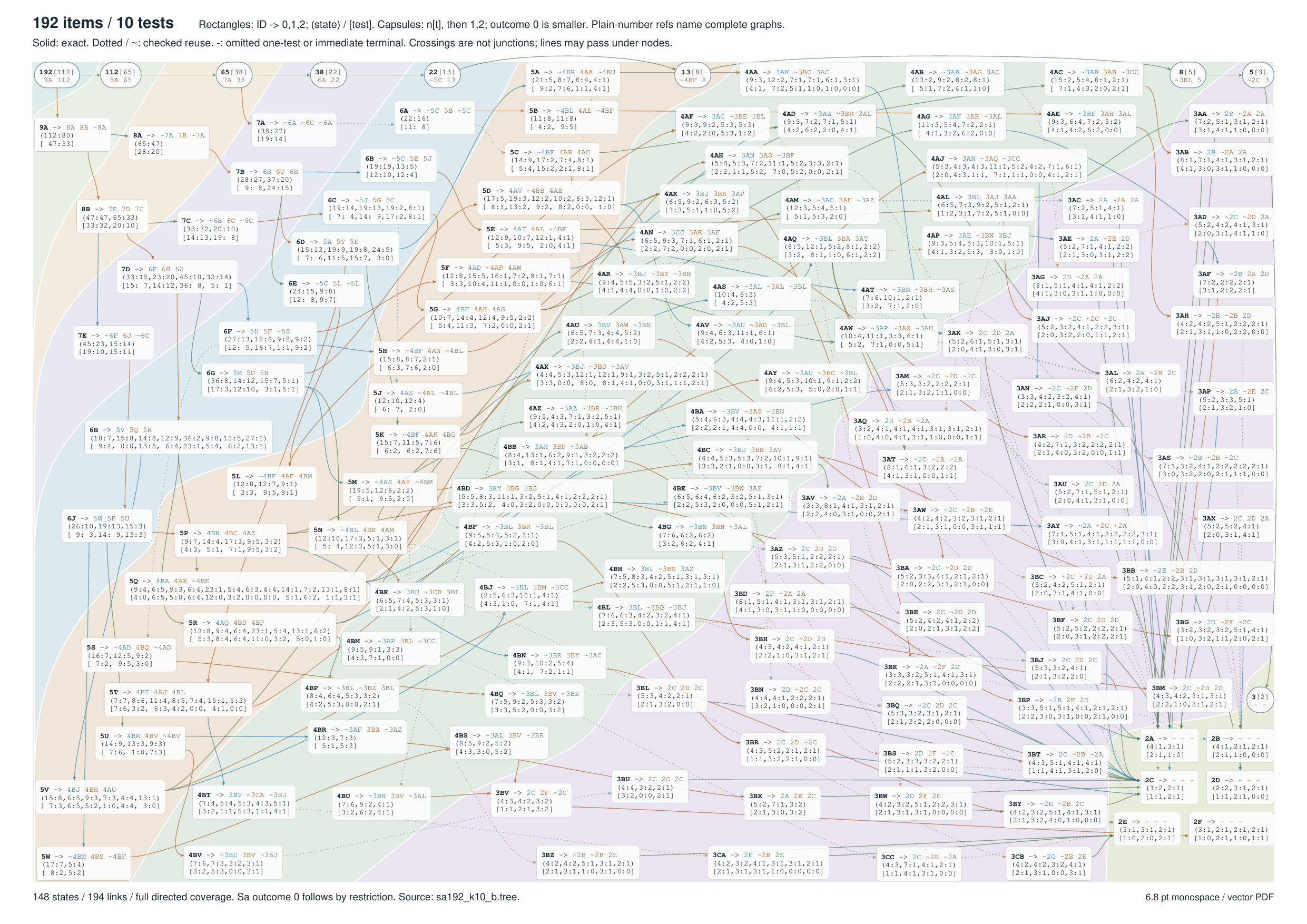}


\begin{thebibliography}{GMSV92}
\small
\bibitem[A86]{A86} M. Aigner, ``Search problems on graphs,'' \emph{Discrete Applied Mathematics} 14 (1986), 215--230. doi: \href{https://doi.org/10.1016/0166-218X(86)90026-0}{\nolinkurl{10.1016/0166-218X(86)90026-0}}.

\bibitem[A88]{A88} M. Aigner, \emph{Combinatorial Search}, Wiley--Teubner, 1988. Section 3.3, especially equation (3.10), display (3.12), and Propositions 3.24--3.25.

\bibitem[GMSV92]{GMSV92} L. Gargano, V. Montuori, G. Setaro and U. Vaccaro, ``An improved algorithm for quantitative group testing,'' \emph{Discrete Applied Mathematics} 36 (1992), 299--306. doi: \href{https://doi.org/10.1016/0166-218X(92)90260-H}{\nolinkurl{10.1016/0166-218X(92)90260-H}}.

\bibitem[ZBM87]{ZBM87} Z. Zhang, T. Berger and J. L. Massey, ``Some Families of Zero-Error Block Codes for the Two-User Binary Adder Channel with Feedback,'' \emph{IEEE Transactions on Information Theory} 33(5) (1987), 613--619. doi: \href{https://doi.org/10.1109/TIT.1987.1057358}{\nolinkurl{10.1109/TIT.1987.1057358}}.

\bibitem[H87]{H87} F. K. Hwang, ``A Tale of Two Coins,'' \emph{American Mathematical Monthly} 94(2) (1987), 121--129. doi: \href{https://doi.org/10.2307/2322412}{\nolinkurl{10.2307/2322412}}.

\bibitem[H89]{H89} F. K. Hwang, ``Updating a Tale of Two Coins,'' \emph{Annals of the New York Academy of Sciences} 576 (1989), 259--265. doi: \href{https://doi.org/10.1111/j.1749-6632.1989.tb16406.x}{\nolinkurl{10.1111/j.1749-6632.1989.tb16406.x}}.

\bibitem[H90]{H90} F. H. Hao, ``The optimal procedures for quantitative group testing,'' \emph{Discrete Applied Mathematics} 26 (1990), 79--86. doi: \href{https://doi.org/10.1016/0166-218X(90)90022-5}{\nolinkurl{10.1016/0166-218X(90)90022-5}}.

\bibitem[C94]{C94} C. A. Christen, ``Search problems: one, two or many rounds,'' \emph{Discrete Mathematics} 136 (1994), 39--51. doi: \href{https://doi.org/10.1016/0012-365X(94)00106-S}{\nolinkurl{10.1016/0012-365X(94)00106-S}}.

\bibitem[LWT18]{LWT18} S. Li, X. Wu and E. Triesch, ``A ternary search problem on two disjoint sets,'' \emph{Discrete Applied Mathematics} 251 (2018), 221--235. doi: \href{https://doi.org/10.1016/j.dam.2018.05.026}{\nolinkurl{10.1016/j.dam.2018.05.026}}.

\bibitem[FHJ22]{FHJ22} S. H. Florin, M. H. Ho and Z. Jiang, ``On the Binary Adder Channel With Complete Feedback, With an Application to Quantitative Group Testing,'' \emph{IEEE Transactions on Information Theory} 68(5) (2022), 2839--2856. doi: \href{https://doi.org/10.1109/TIT.2021.3137965}{\nolinkurl{10.1109/TIT.2021.3137965}}.

\bibitem[Kar26a]{Kar26a} F. Karpelevitch, \emph{Counterexamples to Aigner\textquotesingle s majorization conjecture for star-forest search}, companion manuscript, 2026.
\end{thebibliography}
\end{document}